\documentclass[reqno,A4paper]{amsart}

\usepackage{amsmath}
\usepackage{amssymb}
\usepackage{amsthm}
\usepackage{enumerate}
\usepackage{mathrsfs} 
\usepackage{eqlist}
\usepackage{array}

\theoremstyle{plain}
\newtheorem{theorem}{Theorem}[section]
\newtheorem{prop}[theorem]{Proposition}
\newtheorem{lemma}{Lemma}[section]
\newtheorem{corol}{Corollary}[theorem]

\theoremstyle{definition}

\numberwithin{equation}{section}

\begin{document}

\title[Direct and inverse approximation theorems]%
{Direct and inverse approximation theorems \\ of functions in the Orlicz type spaces ${\mathcal S}_{_{\scriptstyle  M}}$}
\author[Stanislav~Chaichenko, Andrii~Shidlich \and Fahreddin~Abdullayev]%
{Stanislav~Chaichenko*, Andrii~Shidlich** \and Fahreddin~Abdullayev***}

\newcommand{\acr}{\newline\indent}

\address{\llap{*\,}Donbas State Pedagogical University,\acr
                   19, G.~Batyuka st., Slaviansk, \acr
                   Donetsk region, UKRAINE, 84116 }
\email{s.chaichenko@gmail.com}

\address{\llap{**\,}Department of Theory of Functions,\acr
                   Institute of Mathematics of NAS of Ukraine,\acr
                   3, Tereshchenkivska str., Kyiv, 
                   UKRAINE, 01601 }
\email{shidlich@gmail.com}

\address{\llap{***\,}Faculty of Sciences, 
                   Kygyz-Turkish Manas University,\acr
                   56, Chyngyz Aitmatov avenue, Bishkek, 
                   KYRGYZ REPUBLIC, 720044;\acr
                   Faculty of Science and Letters, Mersin University,\acr
                   \c{C}iftlikk\"{o}y Kamp\"{u}s\"{u}, Yeni\c{s}ehir, Mersin, TURKEY, 33342}
\email{fahreddinabdullayev@gmail.com }


\thanks{This work was supported in part by the Ministry of Education and Science of Ukraine within the framework of the fundamental research No.~0118U003390 
and  the Kyrgyz-Turkish Manas University (Bishkek / Kyrgyz Republic), project No.~KTM\"{U}-BAP-2018.FBE.05.}


\begin{abstract}
In the Orlicz type spaces ${\mathcal S}_{M}$, we prove direct and inverse approximation theorems
in terms of the best approximations of functions and moduli of smoothness of fractional order.  We also show
the equivalence between moduli of smoothness and  Peetre $K$-functionals in the spaces ${\mathcal S}_{M}$.
\end{abstract}

\maketitle
\normalsize

\section{Introduction}

Direct approximation theorems are statements asserting that
smoothness of the function $f$ implies a quick decrease to zero of its error of approximation by polynomials or other approximating aggregates. On classes of continuously differentiable functions,  such theorems were first proved in terms of the first-order modulus of continuity  by Jackson  \cite{Jackson_1911}  in 1911.  Later, Zygmund \cite{Zygmund_1924}  and Akhiezer \cite{Akhiezer_M1947} generalized Jackson's results to the second-order modulus of continuity, and Stechkin \cite{Stechkin_1951} extended these results to the moduli of continuity of an arbitrary integer order  $k$, $k\ge 3.$

Inverse approximation theorems are the converse statements that characterize the smoothness
properties of a function depending on the speed of convergence to zero of its approximation by some approximating aggregates. These theorems were first obtained by Bernstein \cite{Bernstein_1912} in 1912. And  already in 1919, direct and inverse  approximation  theorems, due to Jackson and Bernstein, were  given
in the book on approximation theory by  de la Vall$\rm \acute{e}$e Poussin  \cite{Vallee-Poussin_M1918}.

Investigations of the connection (direct and inverse) between the smoothness properties of functions and the possible
orders of their approximations were carried out by many authors on various classes of functions and
for various approximating aggregates. Such results constitute the classics
of modern approximation theory and they are also described quite fully in the monographs \cite{Butzer_Nessel_M1971}, \cite{DeVore_Lorentz_M1993}, \cite{Dzyadyk_Shevchuk_M2008}, \cite{M_Timan_M2009}.

For the last decades, in addition to the classical  direction of theory of direct and  inverse approximation theorems, a number of "non-classical"\ directions have also been developed fruitfully. It should be mention the studies on direct and inverse approximation theorems in the Orlicz function spaces,
which results are contained, in particular, in the papers of Ramazanov \cite{Ramazanov_1984}, Garidi \cite{Garidi_1991}, Runovski \cite{Runovski_2001}, Israfilov and Guven \cite{Israfilov_Guven_2006}, \cite{Guven_Israfilov_2012},   Akg\"{u}n and  Izrafilov \cite{Akgun_Izrafilov_2011},  ~Akg\"{u}n
\cite{Akgun_2012}, Chaichenko \cite{Chaichenko_2013}  and others.

In 2001, Stepanets \cite{Stepanets_2001} considered the spaces  ${\mathcal S}^p$ of $2\pi$-periodic  Lebesgue summable functions $f$ ($f\in L$) with the finite norm
$$
    \|f\|_{_{\scriptstyle {\mathcal S}^p}}=\|\{\widehat f({k})\}_{{k}\in\mathbb  Z}
    \|_{l_p({\mathbb Z})}=\Big(\sum_{{ k}\in\mathbb  Z}|\widehat f({k})|^p\Big)^{1/p},
 $$
where
$\widehat{f}(k):={[f]}\widehat{\ \ }(k)=\frac 1{2\pi}\int_0^{2\pi}f(t) \mathrm{e}^{- \mathrm{i}kt}\mathrm{d}t$,
$k\in\mathbb Z$, are the Fourier coefficients of the function $f$,
and  investigated some approximation characteristics of these spaces,
including in the context of direct and inverse theorems.
Stepanets and Serdyuk \cite{Stepanets_Serdyuk_2002} introduced the notion of $k$th modulus of smoothness in ${\mathcal S}^p$
and established direct and inverse theorems on approximation in terms of these moduli of
smoothness  and the best approximations of  functions. Also this topic was
investigated actively in \cite{Sterlin_1972}, \cite{Vakarchuk_2004}, \cite{Stepanets_2006}, \cite[Ch.~3]{M_Timan_M2009} and others.

In the papers  \cite{Shidlich_Chaichenko_2004}, \cite{Shidlich_Chaichenko_2005} some results for the spaces ${\mathcal S}^p$ were extended  to the Orlicz sequence  spaces.
In particular, in \cite{Shidlich_Chaichenko_2005} the authors found  the exact values of the best $n$-term approximations
and Kolmogorov widths of certain sets of images of the diagonal operators in
the Orlicz spaces. The purpose of this paper is to combine the
above mentioned studies and prove direct and inverse theorems in the Orlicz type
spaces ${\mathcal S}_{M}$  in terms of best approximations of functions and
moduli of smoothness of fractional order.


\section{Preliminaries}

An Orlicz function  $M(t)$ is a non-decreasing convex  function
defined for  $t\ge 0$ such that
$M(0)=0$ and $M(t)\to \infty$  as $t\to \infty$.
Let ${\mathcal S}_{M}$ be the space of all functions $f\in L$ such that
the following quantity
(which is also called the Luxemburg norm of $f$) is finite:
\begin{equation}\label{S_M.1}
    \|{f}\|_{_{\scriptstyle  M}}:=
    \|\{\widehat{f}(k)\}_{k\in {\mathbb Z}}\|_{_{\scriptstyle l_M({\mathbb Z})}}=
    \inf\bigg\{a>0:\  \sum\limits_{k\in\mathbb Z}  M(|{\widehat{f}(k)}|/{a})\le 1\bigg\}.
\end{equation}
Functions $f\in L$ and $g\in L$ are equivalent in the space ${\mathcal S}_{M}$, when $\|f-g\|_{_{\scriptstyle  M}}\!=\!0.$

The spaces ${\mathcal S}_{M}$ defined in this way are Banach spaces.
In case $M(t)=t^p$, $p\ge 1$, they coincide with the above-defined spaces ${\mathcal S}^p$.

Let ${\mathcal T}_{n}$, $n=0,1,\ldots$, be the set of all\ \  trigonometric\ \  polynomials
$\tau_{n}(x){:=}\sum_{|k|\le n}  c_{k}\mathrm{e}^{ikx}$ of the order $n$, where $c_{ k}$
are arbitrary complex numbers. For any function $f\in {\mathcal S}_{M}$, we denote by
\begin{equation}\label{S_M.3}
    E_n (f)_{_{\scriptstyle  M}}:=
    \inf\limits_{\tau_{n-1}\in {\mathcal T}_{n-1} }\|f-\tau_{n-1}\|_{_{\scriptstyle  M}}
\end{equation}
the best approximation of $f$ by trigonometric polynomials $\tau_{n-1}\in {\mathcal T}_{n-1}$ in the space ${\mathcal S}_{M}$.

Similarly to \cite{Butzer_Westphal_1975},  we define the (right) difference of $f\in L$ of fractional order  $\alpha>0$
with respect to the increment $h\in {\mathbb R}$ by
\begin{equation}\label{S_M.6}
    \Delta_h^\alpha f({ x}):=\sum\limits_{j=0}^\infty (-1)^j {\alpha \choose j} f({ x}-jh),\ \
\end{equation}
where ${\alpha \choose j}=\frac {\alpha(\alpha-1)\cdot\ldots\cdot(\alpha-j+1)}{j!},~ j \in \mathbb{N}$, ${\alpha \choose 0}:=1$, and assemble some basic properties of the fractional differences.

{\lemma\label{Lemma_1} Assume that $f\in {\mathcal S}_{M}$, $\alpha, \beta>0$, $x,h\in {\mathbb R}$. Then

\noindent {\rm (i)} $\|\Delta_h^\alpha f\|_{_{\scriptstyle  M}}\le K(\alpha)\|f\|_{_{\scriptstyle  M}}$, \ where $K(\alpha):=\sum_{j=0}^\infty |{\alpha \choose j}|\le 2^{\{\alpha\}}$,
$\{\alpha\}=\inf\{k\in {\mathbb N}: k\ge \alpha\}$.

\noindent {\rm (ii)} ${[\Delta_h^\alpha f]}\widehat {\ \ }(k)=(1-\mathrm{e}^{-ikh})^\alpha \widehat{f}(k)$, $ k\in\mathbb Z$.

\noindent {\rm (iii)} $(\Delta_h^\alpha (\Delta_h^\beta f))(x)=\Delta_h^{\alpha+\beta} f(x)$ (a.\,e.).

\noindent {\rm (iv)} $\|\Delta_h^{\alpha+\beta} f\|_{_{\scriptstyle  M}}\le 2^{\{\beta\}} \|\Delta_h^{\alpha} f\|_{_{\scriptstyle  M}} $.

\noindent {\rm (v)} $\lim\limits_{|h|\to 0}\|\Delta_h^{\alpha} f\|_{_{\scriptstyle  M}}=0$.
}

\vskip 2mm
The proof of Lemma \ref{Lemma_1} and other auxiliary statements of the paper will be given in Section 7.


Based on definition (\ref{S_M.6}), the modulus of smoothness of $f\in {\mathcal S}_{M}$ of the index $\alpha>0$
is defined by
\[
    \omega_\alpha(f,\delta)_{_{\scriptstyle  M}}:=
    \omega_\alpha(f,\delta)_{_{\scriptstyle {\mathcal S}_M}}=
    \sup\limits_{|h|\le \delta}\|\Delta_h^\alpha f\|_{_{\scriptstyle  M}}.
\]
For convenience, we also assume that $\Delta_h^0 f:=f$ and  $\omega_0(f,\delta)_{_{\scriptstyle  M}}:=\|f\|_{_{\scriptstyle  M}}$.
Using the standard arguments, it can be shown that the functions
$\omega_\alpha(f,\delta)_{_{\scriptstyle  M}}$
possess all the basic properties of ordinary moduli of smoothness. Before formulating them,
we give the definition of the $\psi$-derivative of a function.


Let  $\psi=\{\psi_k\}_{k=-\infty}^\infty$ be an arbitrary sequence of complex numbers, $\psi_k\not=0$, $k\in {\mathbb Z}$. If for a given function $f\in L$ with the Fourier series of the form
$S[f](x)=\sum_{k\in {\mathbb Z}}\widehat {f}(k)\mathrm{e}^{\mathrm{i}kx},$ the series
$\sum_{k\in {\mathbb Z}\setminus\{0\}}\widehat {f}(k)\mathrm{e}^{\mathrm{i}kx}/{\psi_k} $
 is the  Fourier series of a certain function $g\in L$, then  $g$ is called
(see, for example, \cite[Ch.~9]{Stepanets_M2005}) $\psi$-derivative of the function $f$ and is denoted as $g:=f^{\psi}$. It is clear that the Fourier coefficients of functions  $f$ and  $f^{\psi }$  are related by equality
\begin{equation}\label{Fourier coeff}
    \widehat  f(k)=\psi_k\widehat  f^{\psi }(k), \ \ k\in {\mathbb Z}\setminus\{0\}
\end{equation}
and $\widehat  f^{\psi }(0)=0$. In case $\psi_k=|k|^{-r}$, $r>0$, $k\in {\mathbb Z}\setminus\{0\}$, we use the notation $f^{\psi}=:f^{(r)}$.


{\lemma\label{Lemma_2}  Assume that  $f, g\in {\mathcal S}_{M}$, $\alpha\ge \beta>0 $ and $\delta,\,\delta_1,\,\delta_2>0$. Then

\noindent {\rm (i)} $\omega_\alpha(f,\delta)_{_{\scriptstyle  M}}$ is a non-negative increasing continuous  function of\ \ $\delta$ on  $(0,\infty)$ such that

$\!\!\!\!$  $\lim\limits_{\delta\to 0+} \omega_\alpha(f,\delta)_{_{\scriptstyle  M}} =0$.

\noindent {\rm (ii)} $\omega_\alpha(f,\delta)_{_{\scriptstyle  M}}\le 2^{\{\alpha-\beta\}} \omega_\beta(f,\delta)_{_{\scriptstyle  M}}$. 

\noindent {\rm (iii)} $\omega_\alpha(f+g,\delta)_{_{\scriptstyle  M}}\le \omega_\alpha(f,\delta)_{_{\scriptstyle  M}}+\omega_\alpha(g,\delta)_{_{\scriptstyle  M}}$.

\noindent {\rm (iv)}  $\omega _1(f,\delta_1+\delta_2)_{_{\scriptstyle  M}}\le
 \omega _1(f,\delta_1)_{_{\scriptstyle  M}}+\omega _1(f,\delta_2)_{_{\scriptstyle  M}}$.

\noindent  {\rm (v)}  $\omega _\alpha(f,\delta)_{_{\scriptstyle  M}}\le 2^{\{\alpha\}}\|f\|_{_{\scriptstyle  M}}$.

\noindent  {\rm (vi) }  if there exists a derivative $f^{(\beta)} \in {\mathcal S}_{M}$, then
  $\omega _{\alpha}(f,\delta)_{_{\scriptstyle  M}}\le \delta^{\beta}\omega _{\alpha-\beta}(f^{(\beta)},\delta)_{_{\scriptstyle  M}}$.

\noindent {\rm (vii) }
 \ $\omega_\alpha(f, p\delta)_{_{\scriptstyle  M}} \le p^\alpha \omega_\alpha(f, \delta)_{_{\scriptstyle  M}}$ \quad $({\alpha\in \mathbb{N}},\ \ {p \in \mathbb{N}})$.

\noindent  {\rm (viii) }
 \ $\omega_\alpha(f, \eta )_{_{\scriptstyle  M}} \le \delta^{-\alpha} (\delta+\eta)^\alpha  \omega_\alpha(f, \delta)_{_{\scriptstyle  M}}$ \quad $({\alpha\in \mathbb{N}}$).
}


\section{Direct approximation theorems}

{\prop\label{Proposition 1} Let $\psi=\{\psi_k\}_{k=-\infty}^\infty$ be an arbitrary
sequence of complex numbers such that
$\psi_k\not=0$ and $\lim\limits_{|k|\to\infty}|\psi_k|=0$.
If for the function $f\in {\mathcal S}_{M}$  there exists a derivative
$f^{\psi}\!\in\! {\mathcal S}_{M}$, then
$$
    E_n(f)_{_{\scriptstyle  M}} \le \varepsilon_n E_n(f^{\psi})_{_{\scriptstyle  M}},
$$
where $\varepsilon_n=\max\limits_{|k|\ge n} |\psi_k|$. }


\begin{proof}
For a fixed $a>0$ and arbitrary numbers  $c_k\in {\mathbb C}$,
$$
    \sum\limits_{|k|\le n-1} M({|\widehat{f}(k)-c_k|}/{a})+ \sum\limits_{|k|\ge n} M({|\widehat{f}(k)|}/{a})\ge \sum\limits_{|k|\ge n} M({|\widehat{f}(k)|}/{a}),
 $$
therefore, for any function $f\in {\mathcal S}_{M}$  we have
\begin{equation}\label{S_M.4}
    E_n (f)_{_{\scriptstyle  M}}=\|f-{S}_{n-1}({f})\|_{_{\scriptstyle  M}}=
    \inf\bigg\{a>0: \sum\limits_{|k|\ge n} M({|\widehat{f}(k)|}/{a})\le 1\bigg\},
\end{equation}
where
$S_{n-1}(f,x)= \sum _{|k|\le n-1}\widehat{f}(k) {\mathrm{e}^{\mathrm{i}kx}}$
is the Fourier sum of the function $f$.


According to  (\ref{S_M.4}) and  (\ref{Fourier coeff}), we have
\begin{eqnarray}\nonumber
 E_n (f)_{_{\scriptstyle  M}}&=&    
    \inf\bigg\{a>0: \sum\limits_{|k|\ge n} M({|\psi_k\widehat  f^{\psi }(k)|}/{a})\le 1\bigg\}
\\ \nonumber
&\le& \inf\bigg\{a>0: \sum\limits_{|k|\ge n} M({\varepsilon_n|\widehat  f^{\psi }(k)|}/{a})\le 1\bigg\}\le \varepsilon_n E_n(f^{\psi})_{_{\scriptstyle  M}}.
\end{eqnarray}
In this case, if  $\varepsilon_n=\max\limits_{|k|\ge n} |\psi_k|=|\psi_{k_0}|$, where $k_0$ is   an integer, $|k_0|\ge n$, then for an arbitrary polynomial of the form $\tilde{\tau}_{k_0}(x):=c\,\mathrm{e}^{{\mathrm i}k_0x}$, $c\not =0$, obviously, the equality holds
$$
    E_n(\tilde{\tau}_{k_0})_{_{\scriptstyle  M}}= \varepsilon_n
    E_n(\tilde{\tau}_{k_0}^{\psi})_{_{\scriptstyle  M}}.
$$
\end{proof}

{\theorem\label{Theorem_1} If $ f\in {\mathcal S}_{M},$  then for any numbers $\alpha>0$ and
 $n \in \mathbb{N}$  the following inequality holds:
\begin{equation}\label{En<omega}
    E_n (f)_{_{\scriptstyle  M}} \le C(\alpha)\, \omega _\alpha(f, {n^{-1}})_{ M}.
\end{equation}
where $C=C(\alpha)$ is a constant that does not depend on $f$ and $n.$}

Before proving Theorem \ref{Theorem_1}, we formulate the auxiliary Lemma \ref{Lemma_3}. This assertion
establishes the equivalence of the Luxemburg norm (\ref{S_M.1}) and the Orlicz norm,
where the latter is defined as follows. Consider the function
\begin{equation} \label{M_tilde}
    \tilde{M}(v):=\sup\{uv-M(u): ~u\ge 0\}
\end{equation}
and the set $\Lambda=\Lambda(\tilde{M})$ of all sequences of positive numbers $\lambda=\{\lambda_k\}_{k\in \mathbb{Z}}$
such that  $\sum_{k\in \mathbb{Z}}\tilde{M}(\lambda_k)\le 1$. For any function  $f\in {\mathcal S}_{M}$, define its Orlicz norm by the equality
\begin{equation} \label{def-Orlicz-norm}
    \|f\|^\ast_{_{\scriptstyle  M}}:= \sup \Big\{ \sum\limits_{k \in \mathbb{Z}}
    \lambda_k|\widehat{f}(k) |: \quad  \lambda\in \Lambda\Big\}.
\end{equation}

{\lemma\label{Lemma_3} For any function $f \in {\mathcal S}_{M}$, the following relation holds:
\begin{equation} \label{estim-for-norms}
    \| f\|_{_{\scriptstyle  M}} \le \| f\|^\ast_{_{\scriptstyle  M}}\le 2 \,\| f\|_{_{\scriptstyle  M}}.
\end{equation}}


\noindent\textit{  Proof of Theorem \ref{Theorem_1}.}  Let us use the proof scheme from  \cite{Stechkin_1951}. Let $\{K_n(t)\}_{n=1}^\infty$ be a sequence of kernels (where $K_n(t)$ is a trigonometric polynomial of order not greater than $n$) such that for all $n=1,2,\ldots$ the following conditions are fulfilled:
\begin{equation} \label{K-cond-1}
    \int_{-\pi}^\pi K_n(t)~\mathrm{d}t=1,
\end{equation}
\begin{equation} \label{K-cond-3}
    \int_{-\pi}^\pi |t|^r |K_n(t)|~\mathrm{d}t\le C(r) (n+1)^{-r} , \quad  r=0,1,2,\ldots.
\end{equation}
In the role of such kernels, in particular, we can take  the well-known Jackson kernels of sufficiently great order, that is,
$$
    K_n(t)=b_p\Big(\frac{\sin pt/2}{\sin t/2}\Big)^{2k_0},
$$
where $k_0$ is an integer that does not depend on  $n,~ 2k_0\ge r+2,$ the positive integer $p$ is
determined from the inequality
${n}/{(2k_0)}<p\le {n}/{(2k_0)}+1, $
and the constant $b_p$ is chosen due to the normalization condition (\ref{K-cond-1}).

It was shown in \cite{Stechkin_1951}  that for any sequence of kernels  $\{K_n(t)\}$  satisfying conditions (\ref{K-cond-1})--(\ref{K-cond-3}), the following estimate holds:
\begin{equation} \label{K-cond-4}
    \int
    _{-\pi}^\pi (|t|+n^{-1})^r ~|K_n(t)|~\mathrm{d}t\le C(r) n^{-r} , \quad (r, n=1,2,\ldots).
\end{equation}

Let us first consider the case of $\alpha \in \mathbb{N}$. Set
$$
    \sigma_{n-1}(x)= (-1)^{\alpha+1} \int
    _{-\pi}^\pi K_{n-1}(t) \sum
    _{j=1}^\alpha (-1)^j {\alpha \choose j} f({ x}-jt)~\mathrm{d}t.
$$
It is clear that $\sigma_{n-1}(x)$ is a trigonometric polynomial
which order does not exceed $n$. Further, in view of  (\ref{K-cond-1}), we have
$$
 f(x)-\sigma_{n-1}(x)= (-1)^{\alpha}\int_{-\pi}^\pi K_{n-1}(t)
    \sum\limits_{j=0}^\alpha (-1)^j {\alpha \choose j} f({ x}-jt)\mathrm{d}t=
    (-1)^{\alpha}\int_{-\pi}^\pi K_{n-1}(t) \Delta_t^\alpha f(x)\mathrm{d}t.
$$
Hence, taking into account  relations (\ref{def-Orlicz-norm})--(\ref{estim-for-norms}) and the definition of the set $\Lambda$, we obtain
\begin{eqnarray}\nonumber
\lefteqn{
    E_n (f)_{_{\scriptstyle  M}} \le \| f-\sigma_{n-1} \|_{_{\scriptstyle  M}} \le
    \| f-\sigma_{n-1} \|^\ast_{_{\scriptstyle  M}}=
     \Big\|(-1)^{\alpha}\int_{-\pi}^\pi K_{n-1}(t) \Delta_t^\alpha f ~\mathrm{d}t  \Big\|^\ast_{_{\scriptstyle  M}}} \\ \nonumber
&=& \sup\Big\{\sum\limits_{k \in \mathbb{Z}} \lambda_k\Big|\frac{1}{2\pi}
    \int_{-\pi} ^\pi  \bigg(\int_{-\pi} ^\pi
    K_{n-1}(t) \Delta_t^\alpha f(x)~ \mathrm{d}t\bigg)~\mathrm{e}^{-\mathrm{i}kx}~{\mathrm d}x  \Big|: \quad
    \lambda \in \Lambda \Big\}.
\end{eqnarray}
Applying now the Fubini theorem and again using estimate (\ref{estim-for-norms}), we find
\begin{eqnarray}\nonumber
    E_n (f)_{_{\scriptstyle  M}} &\le&   \int_{-\pi} ^\pi |K_{n-1}(t)|
    \sup\Big\{\sum\limits_{k \in \mathbb{Z}}
    \lambda_k\Big| \frac{1}{2\pi} \int_{-\pi} ^\pi
    \Delta_t^\alpha f(x)~\mathrm{e}^{-\mathrm{i}kx}~{\mathrm d}x \Big| :~~
    \lambda\in \Lambda \Big\}~\mathrm{d}t
    \\ \nonumber
    &\le&  2 \int_{-\pi} ^\pi|K_{n-1}(t)|~ \|\Delta_t^\alpha f(x)\|_{_{\scriptstyle  M}}^\ast~ \mathrm{d}t \le
    2 \int_{-\pi} ^\pi|K_{n-1}(t)|~ \|\Delta_t^\alpha f(x)\|_{_{\scriptstyle  M}}~ \mathrm{d}t
    \\ \label{int-K-omega}
    &\le&
    2\int_{-\pi} ^\pi |K_{n-1}(t)| \omega_\alpha(f,|t|)_{_{\scriptstyle  M}}~\mathrm{d}t.
\end{eqnarray}
To estimate the integral on the right-hand side of (\ref{int-K-omega}), we use the property (viii) of Lemma \ref{Lemma_2}. Setting $\eta=|t|$, $\delta=n^{-1}$, we see that
 $$
 \omega_\alpha(f; |t|)_{_{\scriptstyle  M}} \le n^\alpha (|t|+n^{-1})^\alpha \omega_\alpha(f,n^{-1})_{_{\scriptstyle  M}}.
 $$
  This inequality together with (\ref{K-cond-4}) yields
  $$
  \int_{-\pi} ^\pi |K_{n-1}(t)| \omega_\alpha(f,|t|)_{_{\scriptstyle  M}}\mathrm{d}t \le
    n^\alpha  \omega_\alpha(f, n^{-1})_{_{\scriptstyle  M}} \int_{-\pi} ^\pi (|t|+n^{-1})^\alpha |K_{n-1}(t)| \mathrm{d}t\le C(\alpha) \omega_\alpha (f, {n^{-1}} )_{ M}.
  $$
Thus, in the case of $\alpha \in \mathbb{N}$, the theorem is proved.

If $\alpha >0$, $\alpha \not \in \mathbb{N},$ then we denote by $\beta$ an arbitrary positive integer satisfying the condition $\beta-1<\alpha<\beta$. Due to property (ii) of Lemma \ref{Lemma_2}, we obtain
$$
    E_n (f)_{_{\scriptstyle  M}} \le C(\beta)~ \omega_\beta (f, {n^{-1}} )_{ M} \le
    C(\beta)~ \omega_\alpha (f, {n^{-1}} )_{ M}.
$$
\vskip -3mm$\hfill\Box$


\section{ Inverse approximation theorems}  The key role in
proving of the inverse approximation  theorems is played by the known Bernstein  inequality in
which the norm of the derivative of a trigonometric polynomial is estimated  in terms of the norm
of this polynomial (see, e.g. \cite[Ch.~4]{A_Timan_M1960}, \cite[Ch.~4]{M_Timan_M2009}).


{\prop\label{Proposition 2} Let $\psi=\{\psi_k\}_{k=-\infty}^\infty$ be an arbitrary sequence of complex numbers, $\psi_k \not=0$. Then for any $\tau_n\in {\mathcal T}_{n}$, $n\in \mathbb{N}$, the following inequality holds:
\begin{equation}\label{Bernstein_Ineq}
    \|\tau^\psi_n \|_{_{\scriptstyle  M}}\le \frac 1{\epsilon_n}\|\tau_n\|_{_{\scriptstyle  M}}, \quad
    \epsilon_n:=\min_{0<| k | \le n}|\psi_k|,
\end{equation}}

\begin{proof}
 Let $\tau_n(x)=\sum_{|k|\le n}  c_{k}\mathrm{e}^{\mathrm{i}kx}$, $c_k\in {\mathbb C}$.  By the definition of the $\psi$-derivative and equalities (\ref{Fourier coeff}), we get
\begin{eqnarray}\nonumber
\lefteqn{    \|\tau^\psi_n \|_{_{\scriptstyle  M}}= 
       \inf \Big\{a>0:~  \sum\limits_{0<| k | \le n} M\Big({|c_k| /
       |a\psi_k |}\Big)\le 1\Big\}}  \\ \nonumber
       &\le&  \max_{0<| k | \le n}{|\psi_k|^{-1}} \inf \Big\{a>0: ~ \sum\limits_{0<| k | \le n}
      M\Big({| c_k| /  a}\Big)\le 1 \Big\}= \frac 1{\epsilon_n} \|\tau_n\|_{_{\scriptstyle  M}}.
\end{eqnarray}
\end{proof}
In this case, if  $\epsilon_n=\min\limits_{0<| k | \le n} |\psi_k|=|\psi_{k_0}|$,   then for an arbitrary polynomial of the form $\tilde{\tau}_{k_0}(x):=c\,\mathrm{e}^{\mathrm{i}k_0x}$, $c\not =0$, we have
$$
    \|\tilde{\tau}_{k_0}^\psi \|_{_{\scriptstyle  M}}=
    \inf \Big\{a>0: \ M\Big({|c|}/{|a\psi_{k_0}|}\Big)\le 1 \Big\}=
    \frac 1{\epsilon_n} \| \tau_{k_0} \|_{_{\scriptstyle  M}}.
$$


{\corol\label{Corollary 1} Let $\psi=\{\psi_k\}_{k=-\infty}^\infty$ be an arbitrary sequence of complex numbers {such that}  $|\psi_{-k}|=|\psi_k|\ge|\psi_{k+1}|>0$.  Then for any $\tau_n\in {\mathcal T}_{n}$, $n\in \mathbb{N}$,
\begin{equation}\label{Bernstein_Ineq1}
    \|\tau^\psi_n \|_{_{\scriptstyle  M}}\le \frac 1{|\psi_n|}\|\tau_n\|_{_{\scriptstyle  M}}.
\end{equation}
In particular, if $\psi_k=|k|^{-r}$, $r>0$, $k\in {\mathbb Z}\setminus\{0\}$, then
$$
    \|\tau^\psi_n \|_{_{\scriptstyle  M}} = \|\tau^{(r)}_n \|_{_{\scriptstyle  M}} \le
    n^r \|\tau_n\|_{_{\scriptstyle  M}}.
$$}

{\theorem\label{Theorem_2} If $ f\in {\mathcal S}_{M}$, then for any $\alpha>0$ and $n\in {\mathbb N}$, the following inequality is true:
\begin{equation}\label{S_M.12}
    \omega _\alpha(f, {n^{-1}} )_{\scriptstyle  M}\le \frac {C(\alpha) }{n^\alpha}
    \sum _{\nu =1}^n \nu ^{\alpha-1} E_{\nu}(f)_{_{\scriptstyle  M}},
\end{equation}
where $C=C(\alpha)$ is a constant that does not depend on $f$ and $n.$
}

For the spaces $L_p$ of $2\pi$-periodic functions integrable to the $p$th power with the usual norm, inequalities of the type (\ref{S_M.12}) were proved in \cite{M_Timan_1958} (see also \cite[Ch.~6]{A_Timan_M1960}, \cite[Ch.~2]{M_Timan_M2009}). In the spaces ${\mathcal S}^p$, similar results were obtained in \cite{Stepanets_Serdyuk_2002}, \cite{Sterlin_1972}.

\begin{proof}
 Let us use the proof scheme from \cite[Ch.~6]{A_Timan_M1960}. Let $ f\in {\mathcal S}_{M}$ and ${S}_n:={S}_n(f)$  be the Fourier sum of $f$.  Then, due to Lemma \ref{Lemma_2} (v)  and relation (\ref{S_M.4}) for an arbitrary $m\in {\mathbb N}$, we have
\begin{eqnarray}\nonumber
    \omega _\alpha(f, {n^{-1}} )_{ M} &\le&
    \omega _\alpha(f-S_{2^{m+1}}, {n^{-1}})_{ M}+
    \omega _\alpha(S_{2^{m+1}}, {n^{-1}})_{ M}
      \\ \label{S_M.9}
    &\le&
    2^{\{\alpha\}} E_{2^{m+1}+1}(f)_{_{\scriptstyle  M}}+
    \omega _\alpha(S_{2^{m+1}}, {n^{-1}})_{ M}.
\end{eqnarray}
Further, using  property {\rm (vi)} of Lemma \ref{Lemma_2} and the properties of the norm, we obtain
$$
    \omega _\alpha(S_{2^{m+1}}, {n^{-1}})_{ M}\le n^{-\alpha}\|S_{2^{m+1}}^{(\alpha)}\|_{_{\scriptstyle  M}}\le
    n^{-\alpha} \Big(\|S_{1}^{(\alpha)}\|_{_{\scriptstyle  M}}+
    \sum_{k=0}^{m} \|S_{2^{k+1}}^{(\alpha)}-S_{2^{k}}^{(\alpha)}\|_{_{\scriptstyle  M}}\Big).
$$
Moreover, on the basis of Corollary \ref{Corollary 1},
$$
    \|S_{2^{k+1}}^{(\alpha)}-S_{2^{k}}^{(\alpha)}\|_{_{\scriptstyle  M}}\le {2^{(k+1)\alpha}}\|S_{2^{k+1}}-S_{2^{k}}\|_{_{\scriptstyle  M}}\le
    2^{(k+1)\alpha+1}E_{2^{k}+1}(f)_{_{\scriptstyle  M}},
$$
and
$\|S_{1}^{(\alpha)}\|_{_{\scriptstyle  M}}=\|S_{1}^{(\alpha)}-S_0^{(\alpha)}\|_{_{\scriptstyle  M}}
\le 2E_{1}(f)_{_{\scriptstyle  M}}.$
Therefore,
$$
    \omega _\alpha(S_{2^{m+1}}, {n^{-1}})_{ M}\le
    n^{-\alpha}\|S_{2^{m+1}}^{(\alpha)}\|_{_{\scriptstyle  M}}\le
    2n^{-\alpha} \Big(E_{1}(f)_{_{\scriptstyle  M}}+
    \sum_{k=0}^{m}  2^{(k+1)\alpha}E_{2^{k}+1}(f)_{_{\scriptstyle  M}}\Big).
 $$
Taking into account the relation
$$
  2^{(k+1)\alpha} E_{2^{k}+1}(f)_{_{\scriptstyle  M}} \le 2^{2\alpha}
  \sum_{\nu=2^{k-1}+2}^{2^k+1}\nu ^{\alpha-1}  E_{\nu}(f)_{_{\scriptstyle  M}},\quad k=1,2,\ldots,
$$
we get
\begin{eqnarray}\nonumber
    \omega _\alpha(S_{2^{m+1}}, {n^{-1}})_{ M}&\le&
    2^{2\alpha+1} n^{-\alpha}\Big(E_{1}(f)_{_{\scriptstyle  M}}+E_{2}(f)_{_{\scriptstyle  M}}+
    \sum_{k=1}^{m} \sum_{\nu=2^{k-1}+2}^{2^k+1} \nu ^{\alpha-1}  E_{\nu}(f)_{_{\scriptstyle  M}}\Big)
    \\ \nonumber
    &\le&
    \frac{c^*(\alpha)}{n^{\alpha}} \sum_{k=1}^{2^m+1} k ^{\alpha-1}  E_{k}(f)_{_{\scriptstyle  M}} .
\end{eqnarray}
Choosing now an integer  $m$ so that $2^m+1\le n\le 2^{m+1}$ and substituting this estimate into (\ref{S_M.9}), we get (\ref{S_M.12}).
\end{proof}

{\corol\label{Corollary_2}  Assume that the sequence of the best approximations $E_n(f)_{_{\scriptstyle  M}}$ of a function $f\in {\mathcal S}_{M}$ satisfies the following relation for some $\beta >0$:
$$
    E_n(f)_{_{\scriptstyle M}}= {\mathcal O}(n^{-\beta }).
$$
Then, for any $\alpha>0$, one has
$$
    \omega _\alpha(f, t)_{_{\scriptstyle M}}=
    \left \{ \begin{matrix}  {\mathcal O}(t^{\beta }) & \hfill \mbox {for} \ \ \beta <\alpha, \hfill \\
                        {\mathcal O}(t^\alpha|\ln t|) & \hfill \mbox {for}\ \ \beta =\alpha, \hfill \\
                        {\mathcal O}(t^\alpha) & \hfill \mbox {for} \ \ \beta >\alpha.\hfill \end{matrix} \right.
$$
}


\section{ Constructive characteristics of the classes of functions defined by the
$\alpha$th moduli of smoothness}

In the following two sections some applications of the obtained results are considered.
In particular, in this section we give the constructive characteristics of the classes
${\mathcal S}_{M}H^{\omega}_{\alpha} $ of functions for which the $\alpha$th moduli of smoothness do not exceed some majorant.

Let $\omega$ be  a function defined on interval $[0,1]$. For a fixed $\alpha>0$, we set
\begin{equation} \label{omega-class}
    {\mathcal S}_{M}H^{\omega}_{\alpha}=
    \Big\{f\in {\mathcal S}_{M}:  \quad \omega_\alpha(f, \delta)_{_{\scriptstyle M}}=
    {\mathcal O}  (\omega(\delta)),\quad  \delta\to 0+\Big\}.
\end{equation}
Further, we consider the functions $\omega(\delta)$, $\delta\in [0,1]$, satisfying the following conditions 1)--4): \textbf{1)} $\omega(\delta)$ is continuous on $[0,1]$;\  \textbf{  2)} $\omega(\delta)\uparrow$;\  \textbf{  3)} $\omega(\delta)\not=0$ for any $\delta\in (0,1]$;\  \textbf{  4)}~$\omega(\delta)\to 0$ as $\delta\to 0+$; \noindent as well-known condition $({\mathcal B}_\alpha)$, $\alpha>0$:
${\sum_{v=1}^n v^{\alpha-1}\omega({v^{-1}}) =
{\mathcal O}  (n^\alpha \omega ( {n^{-1}}))}$ (see, e.g. \cite{Bari_Stechkin_1956}).

{\theorem\label{Theorem_3}  Assume that $\alpha>0$  and $\omega$ is a function, satisfying  conditions  $1)$--\,$4)$ and $({\mathcal  B}_\alpha)$. Then, in order a function $f\in {\mathcal S}_{M}$ to belong to the class ${\mathcal S}_{M}H^{\omega}_{\alpha}$, it is necessary and sufficient that
\begin{equation} \label{iff-theorem}
    E_n(f)_{\scriptstyle  M}={\mathcal O} ( \omega ({n^{-1}} ) ).
\end{equation}
}
\begin{proof}
 Let $f \in {\mathcal S}_{M}H^{\omega}_{\alpha}$, by virtue of Theorem \ref{Theorem_1}, we have
\begin{equation} \label{using-direct-theorem}
    E_n(f)_{\scriptstyle  M} \le C(\alpha) \omega_\alpha (f; {n^{-1}} )_{_{\scriptstyle M}},
\end{equation}
Therefore, {relation (\ref{omega-class}) yields (\ref{iff-theorem})}.
On the other hand, if relation (\ref{iff-theorem}) holds, then by virtue of Theorem \ref{Theorem_2},  taking into account the condition $({\mathcal  B}_\alpha)$, we obtain
\begin{equation} \label{using-inverse-theorem}
    \omega _\alpha(f, {n^{-1}} )_{\scriptstyle  M}\le
    \frac {C(\alpha) }{n^\alpha} \sum _{\nu =1}^n \nu ^{\alpha-1} E_{\nu}(f)_{_{\scriptstyle  M}}
        \le
    \frac {C_1}{n^\alpha} \sum _{\nu =1}^n \nu ^{\alpha-1} \omega ({v^{-1}})=
    {\mathcal O}  (\omega ( {n^{-1}})).
\end{equation}
Thus, the function  $f$ belongs to the set ${\mathcal S}_{M}H^{\omega}_{\alpha}$.
\end{proof}

The function $\varphi (t)=t^r$, $r \le \alpha$, satisfies the condition $({\mathcal  B}_\alpha)$. Hence, denoting by ${\mathcal S}_{M}H_{\alpha}^r$  the class ${\mathcal S}_{M}H^{\omega}_{\alpha}$ for $\omega(t)=t^r$, $0<r\le \alpha,$  we establish the following statement:

\vskip 2mm

{\corol\label{Corollary_5}  Let $\alpha >0$, $0<r\le \alpha.$
In order  a function $f \in S_{\scriptstyle  M}$ to belong to ${\mathcal S}_{M}H_{\alpha}^r$, it is
necessary and sufficient that
$$
    E_n(f)_{\scriptstyle  M}={\mathcal O}   ({n^{-r}} ).
$$
}

\section{The equivalence between $\alpha$th moduli of smoothness and $K$-functionals}

 $K$-functionals were introduced by Lions and Peetre in 1961, and defined in their usual form by Peetre in the monograph \cite{Peetre_M1963} as a basis for his theory of operator interpolation. Unlike the moduli of continuity expressing the smooth properties of functions, $K$-functionals express some of their approximative properties. In this section we prove the equivalence between our moduli of smoothness and certain Peetre $K$-functionals. This connection is important for studying the properties of the modulus of smoothness and the $K$-functional, and also for their further application to the problems of approximation theory.

In the space ${\mathcal S}_{M}$, the Petree $K$-functional of a  function $f$
(see, {e.g.} \cite[Ch.~6]{DeVore_Lorentz_M1993}), {which} generated by its derivative of order $\alpha>0$, is the following quantity:
$$
    K_\alpha(\delta,f)_{_{\scriptstyle  M}}=\inf\Big\{\|f-h\|_{_{\scriptstyle  M}}+
    \delta^\alpha \|h^{(\alpha)}\|_{_{\scriptstyle  M}}:\
    h^{(\alpha)}\in {\mathcal S}_{M}\Big\},\quad \delta>0.
$$

{\theorem\label{Theorem_4} For each $ f\!\in {\mathcal S}_{M}$, $\alpha\!>0$, there exist constants $C_1(\alpha)$, $C_2(\alpha)\!>\!0$, such that for $\delta>0$
\begin{equation}\label{KO1}
      C_1(\alpha)\omega _\alpha(f, \delta)_{_{\scriptstyle  M}}\le
      K_\alpha(\delta,f)_{_{\scriptstyle  M}}\le
      C_2(\alpha)\omega _\alpha(f, \delta)_{_{\scriptstyle  M}}.
\end{equation} }

\begin{proof}
 Consider an arbitrary function $h\in {\mathcal S}_{M}$ such that
$h^{(\alpha)}\in {\mathcal S}_{M}$. Then  we have by
Lemma \ref{Lemma_2} (iii), (v) and (vi)
$$
    \omega _\alpha(f, \delta)_{_{\scriptstyle  M}}\le
    \omega _\alpha(f-h, \delta)_{_{\scriptstyle  M}}+
    \omega _\alpha(h, \delta)_{_{\scriptstyle  M}}\le 2^{\{\alpha\}}\|f-h\|_{_{\scriptstyle  M}}+ \delta^\alpha\|h^{(\alpha)}\|_{_{\scriptstyle  M}}.
$$
Taking the infimum over all  $h\in {\mathcal S}_{M}$ such that
$h^{(\alpha)}\in {\mathcal S}_{M}$,  we get the left-hand side of (\ref{KO1}).

To prove the right-hand side of  (\ref{KO1}),  let us formulate the following auxiliary  lemma.

{\lemma\label{Lemma_4}  Assume that  $\alpha >0$, $n \in \mathbb{N}$ and $0\le h\le 2\pi/n$.
 Then for any  $\tau_n \in {\mathcal T}_n $
\begin{equation}\label{Bermstain-inequl-gener}
    \Big(\frac{ \sin(nh/2)}{n/2} \Big)^\alpha\| \tau_n^{(\alpha)} \|_{_{\scriptstyle  M}}\le
     \|\Delta_h^\alpha \tau_n \|_{_{\scriptstyle  M}}\le h^\alpha\| \tau_n^{(\alpha)} \|_{_{\scriptstyle  M}}.
\end{equation}}

Now let  $\delta \in (0,2\pi)$ and  $n \in \mathbb{N}$ such that
$\pi/n<\delta < 2\pi/n$. Let also $S_n:=S_n(f)$ be the Fourier sum of $f$.    Using Lemma \ref{Lemma_4} with
$h=\pi /n$ and  property  (i) of Lemma \ref{Lemma_1}, we obtain
$$
    \|S_n^{(\alpha)}\|_{_{\scriptstyle  M}} \le
    2^{-\alpha+1} n^{\alpha} \|\Delta_{\pi/n}^\alpha S_n \|_{_{\scriptstyle  M}} 
        \le 2(\pi/\delta)^{\alpha} \Big( \|\Delta_{\pi/n}^\alpha (S_n -f)  \|_{_{\scriptstyle  M}}+ \|\Delta_{\pi/n}^\alpha f  \|_{_{\scriptstyle  M}}\Big) 
$$
\begin{equation}\label{KO2}
    \le 2(\pi/\delta)^{\alpha} \Big( 2^{\{\alpha\}} \| f- S_n \|_{_{\scriptstyle  M}}+
    \|\Delta_{\pi/n}^\alpha f  \|_{_{\scriptstyle  M}}\Big).
\end{equation}
By virtue of (\ref{S_M.4}) and Theorem \ref{Theorem_1}, we have
\begin{equation}\label{KO3}
    \|f- S_n \|_{_{\scriptstyle  M}}=E_n(f)_{_{\scriptstyle  M}} \le
    C(\alpha) \omega_\alpha (f, \delta )_{_{\scriptstyle  M}}.
\end{equation}
Combining (\ref{KO2}), (\ref{KO3})    and the definition of modulus of smoothness, we obtain the relation
$$
    \|S_n^{(\alpha)} \|_{_{\scriptstyle  M}} \le
    C_2(\alpha) \delta^{-\alpha} \omega_\alpha (f, \delta )_{_{\scriptstyle  M}},
$$
where $C_2(\alpha):=2\pi^\alpha(2^{\{\alpha\}}C(\alpha)+1)$, which yields the right-hand side of (\ref{KO1}):
$$
    K_\alpha(\delta,f)_{_{\scriptstyle  M}} \le \|f- S_n \|_{_{\scriptstyle  M}}+
    \delta^\alpha \|S_n^{(\alpha)}\|_{_{\scriptstyle  M}} \le
    C_2(\alpha) \omega_\alpha (f, \delta )_{_{\scriptstyle  M}}.
$$
\end{proof}

\section{Proof of the  auxiliary statements}

\textit{ Proof of Lemma \ref{Lemma_1}.} Let us set $f({x}-jh)=:f_{jh}({x})$. For any  ${ k}\in {\mathbb Z}$  and $j=0,1,\ldots$, we have
 $ \widehat{f}_{jh}({ k})=\widehat{f}({ k})\mathrm{e}^{-{\mathrm i}kjh}$. Therefore,
\begin{eqnarray}\nonumber
\lefteqn{
    \|\Delta_h^\alpha f\|_{_{\scriptstyle  M}}=
    \inf\bigg\{a>0: \sum_{{ k}\in {\mathbb Z}}
    M\Big(|{[\Delta_h^\alpha f]}\widehat {\ \ }(k)|/a\Big)\le 1\bigg\}} \\ \nonumber
&=&\inf\bigg\{a>0: \sum_{{ k}\in {\mathbb Z}} M\Big(\Big|\Big[\sum\limits_{j=0}^\infty
    (-1)^j {\alpha \choose j} f_{jh}\Big]\widehat{\ \ }(k)\Big|/a\Big)\le 1\bigg\}
     \\ \nonumber
&=&\inf\bigg\{a>0: \sum_{{ k}\in {\mathbb Z}} M\Big(\Big|\widehat{f}({ k})\sum\limits_{j=0}^\infty (-1)^j {\alpha \choose j}\mathrm{e}^{-{\mathrm i}kjh}\Big|/a\Big)\le 1\bigg\}.
\end{eqnarray}
For a fixed $a>0$
\begin{eqnarray}\nonumber
    \sum_{{ k}\in {\mathbb Z}} M\Big(\Big|\widehat{f}({ k})\sum\limits_{j=0}^\infty
    (-1)^j {\alpha \choose j}\mathrm{e}^{-{\mathrm i}kjh}\Big|/a\Big)&\le&
    \sum_{{ k}\in {\mathbb Z}}  M\Big(\sum\limits_{j=0}^\infty \Big| {\alpha \choose j}\Big|\,
    |\widehat{f}({ k})|/a\Big)\\ \label{S_M.8}
    &\le&\sum_{{ k}\in {\mathbb Z}}
    M\Big( 2^{\{\alpha\}}  |\widehat{f}({ k})|/a\Big).
\end{eqnarray}
and hence {\rm (i)} holds. Property {\rm (ii)} is obvious:
\begin{equation}\label{difference_Fourier_Coeff}
    {[\Delta_h^\alpha f]}\widehat {\ \ }(k)=\Big[\sum\limits_{j=0}^\infty (-1)^j
    {\alpha \choose j} f_{jh}\Big]\widehat{\ \ }(k)
        =
    \widehat{f}({ k})\sum\limits_{j=0}^\infty
    (-1)^j {\alpha \choose j}\mathrm{e}^{-{\mathrm i}kjh}=(1-\mathrm{e}^{-{\mathrm i}kh})^\alpha \widehat{f}(k),
\end{equation}
and property (iii) is its consequence.  Part {\rm (iv)} follows by
{\rm (i)}--{\rm (iii)}.

To prove {\rm (v)} we first show that the following relation holds:
 \begin{equation}\label{Trigon_Polynomial}
 \lim\limits_{|h|\to 0}\|\Delta_h^{\alpha} \tau_{n}\|_{_{\scriptstyle  M}}=0
\end{equation}
where $\tau_{n}$ is an arbitrary polynomial of the form $\tau_{n}(x){=}\sum_{|k|\le n}  c_{k}\mathrm{e}^{\mathrm{i}kx}$, $n\in {\mathbb N}$, $c_{k}\in {\mathbb C}$.

Since $\|\tau_{n}\|_{_{\scriptstyle  M}}=\inf\{a>0: \sum_{|k|\le n}M( |c_{k}|/a)\le 1\}$, then taking into the account  {\rm (ii)}, for $a_0=|nh|^\alpha \|\tau_{n}\|_{_{\scriptstyle  M}}$, we obtain
\begin{eqnarray}\nonumber
\lefteqn{
    \sum_{|k|\le n} M\Big(|{[\Delta_h^\alpha \tau_{n}]}\widehat {\ \ }(k)|/a_0\Big)=
    \sum_{|k|\le n} M\Big(|1-\mathrm{e}^{-{\mathrm i}kh}|^\alpha |c_{k}|/a_0\Big)  =
    \sum_{|k|\le n} M\Big(2^\alpha \Big|\sin \frac {kh}2\Big|^\alpha |c_{k}|/a_0\Big)}
    \\ \label{Trigon_P}
 &\le&
     \sum_{|k|\le n} M\Big( |kh|^\alpha |c_{k}|/a_0\Big) \le
     \sum_{|k|\le n} M\Big( |nh|^\alpha |c_{k}|/a_0\Big)= \sum_{|k|\le n}M\Big( |c_{k}|/\|\tau_{n}\|_{_{\scriptstyle  M}}\Big)\le 1 .
\end{eqnarray}
Therefore, $\|\Delta_h^{\alpha} \tau_{n}\|_{_{\scriptstyle  M}}\le |nh|^\alpha \|\tau_{n}\|_{_{\scriptstyle  M}}$. For an arbitrary $\varepsilon>0$,  we set $\delta:=\delta(\varepsilon)=
\Big(\varepsilon/n^\alpha\|\tau_{n}\|_{_{\scriptstyle  M}}\Big)^{1/\alpha}$. Then for all $|h|<\delta$, we have
$\|\Delta_h^{\alpha} \tau_{n}\|_{_{\scriptstyle  M}}<\varepsilon$, i.e., relation (\ref{Trigon_Polynomial}) is indeed fulfilled.

Now let $f$ is an arbitrary function from ${\mathcal S}_{M}$ and $S_{n}(f,x)= \sum _{|k|\le n}\widehat{f}(k) {\mathrm{e}^{\mathrm{i}kx}}$ is its Fourier sum. Since the value $\|f-{S}_{n}({f})\|_{_{\scriptstyle  M}}$ tends to zero as $n\to \infty$, then for any  $\varepsilon>0$ there exist a positive integer $n_0=n_0(\varepsilon)$ such that
 for any  $n>n_0$, we have
 \[
 \|f-{S}_{n}({f})\|_{_{\scriptstyle  M}}<{\varepsilon}/ 2^{\{\alpha\}+1},
 \]
 Furthermore, by virtue of  (\ref{Trigon_Polynomial}),   there exist a number  $\delta:=\delta(\varepsilon,n)$ such that $\|\Delta_h^{\alpha} S_{n}(f)\|_{_{\scriptstyle  M}}<\frac {\varepsilon}2$ when $|h|<\delta$. Then using properties of norm and {\rm (i)}, for  $n>n_0$ we get
 \[
 \|\Delta_h^\alpha f\|_{_{\scriptstyle  M}}
 \le\|\Delta_h^\alpha (f-S_n(f))\|_{_{\scriptstyle  M}}+\|\Delta_h^\alpha S_n(f)\|_{_{\scriptstyle  M}}
  \le 2^{\{\alpha\}}\| f-S_n(f)\|_{_{\scriptstyle  M}}+\|\Delta_h^\alpha S_n(f)\|_{_{\scriptstyle  M}}<\varepsilon,
 \]
which yields {\rm (v)}.
$\hfill\Box$

\vskip 2mm \textit{ Proof of Lemma \ref{Lemma_2}.} In {\rm (i)}, the convergence to zero for $\delta\to 0+$ follows by  Lemma \ref{Lemma_1} {\rm (v)}.  Part  {\rm (v)} is the consequence of Lemma \ref{Lemma_1}  {\rm (i)}. Property {\rm (iii)}, non-negativity and increasing of the function $\omega_{\alpha}(f,t)_{_{\scriptstyle  M}}$  follow from the definition of modulus of smoothness. According to Lemma \ref{Lemma_1} {\rm (i)} and {\rm (iii)}, for arbitrary numbers $0<\beta\le \alpha$, we have
 $$
 \|\Delta_h^\alpha f \|_{_{\scriptstyle  M}}=\|\Delta_h^{\alpha-\beta} (\Delta_h^\beta f) \|_{_{\scriptstyle  M}} \le    2^{\{\alpha-\beta\}}\|\Delta_h^\beta f \|_{_{\scriptstyle  M}},
 $$
whence passing to the exact upper bound over all $|h|\le \delta$, we obtain {\rm (ii)}.
Part {\rm (iv)} is proved by  the following standard arguments:
$$
    \omega_1(f,\delta_1+\delta_2)_{_{\scriptstyle  M}}=\!\!\sup\limits_{|h_1|\le \delta_1,|h_2|\le \delta_2} \|f({ x}+h_1+h_2)-f({ x})\|_{_{\scriptstyle  M}}\le
    \sup\limits_{|h_2|\le \delta_2} \|f({ x}+h_1+h_2)-f({ x}+h_1)\|_{_{\scriptstyle  M}}
    $$
    $$
        +\sup\limits_{|h_1|\le \delta_1} \|f({ x}+h_1)-f({ x})\|_{_{\scriptstyle  M}}
        \le  \omega_1(f,\delta_2)_{_{\scriptstyle  M}}+ \omega_1(f,\delta_1)_{_{\scriptstyle  M}}.
$$
In particular, this   yields the continuity of the function $\omega_1(f,\delta)_{_{\scriptstyle  M}}$, since for arbitrary $\delta_1>\delta_2>0$,
$\omega _1(f,\delta_1)_{_{\scriptstyle  M}}-\omega _1(f,\delta_2)_{_{\scriptstyle  M}}\le
\omega_1(\delta_1-\delta_2)_{_{\scriptstyle  M}}\to 0$ as $\delta_1-\delta_2\to 0.$

Let us prove the continuity of 
$\omega_\alpha(f,\delta)_{_{\scriptstyle  M}}$ for arbitrary $\alpha>0$. Let $0\,{<}\delta_1\,{<}\delta_2$  and  $h=h_1+h_2,$ where\ \  $0<h_1\le \delta_1$, $0<h_2\le \delta_2-\delta_1.$ Since
$
    \Delta_h^\alpha f(\delta) =\Delta_{h_1}^\alpha f (\delta) +\sum_{j=0}^\infty {\alpha \choose j}
    (-1)^{j} \Delta_{jh_2}^1 f(\delta-jh_1)
$
and
$$
    \Big\|\sum_{j=0}^\infty {\alpha \choose j}
    (-1)^{j} \Delta_{jh_2}^1 \, f_{jh_1}\Big\|_{_{\scriptstyle  M}}
        =
    \inf\bigg\{a>0: \sum_{{ k}\in {\mathbb Z}}
    M\Big(\Big|{\Big[\sum_{j=0}^\infty {\alpha \choose j}
    (-1)^{j} \Delta_{jh_2}^1  f_{jh_1}\Big]}\widehat {\ \ }(k)\Big|/a\Big)\le 1\bigg\}
$$
$$
    \le \inf\bigg\{a>0: \sum_{{ k}\in {\mathbb Z}}
    M\Big( |2^{\{\alpha\}}\alpha  [\Delta_{h_2}^1 f  ]\widehat {\ \ }(k) |/a\Big)\le 1\bigg\}\le
    2^{\{\alpha\}}\alpha\|\Delta_{h_2}^1 f\|_{_{\scriptstyle  M}},
$$
then
$
    \|\Delta_h^\alpha f\|_{_{\scriptstyle  M}}  \le
    \|\Delta_{h_1}^\alpha f\|_{_{\scriptstyle  M}} +
    2^{\{\alpha\}}\alpha\|\Delta_{h_2}^1 f\|_{_{\scriptstyle  M}}
$
and
$
    \omega_\alpha (f,\delta_2)_{_{\scriptstyle  M}}   \le
    \omega_\alpha (f,\delta_1)_{_{\scriptstyle  M}} +
    2^{\{\alpha\}}\alpha \,\omega_1 (f, \delta_2-\delta_1)_{_{\scriptstyle  M}}.
$
Hence, we obtain the necessary relation:
$$
    \omega_\alpha(f,\delta_2)_{_{\scriptstyle  M}}
    \!- \omega_\alpha(f,\delta_1)_{_{\scriptstyle  M}}\le 2^{\{\alpha\}}\alpha \,
    \omega_1(f, \delta_2-\delta_1)_{_{\scriptstyle  M}} \to 0,\ \ \   \delta_2-\delta_1\to 0.
$$

If there exists a derivative $f^{(\beta)}\in {\mathcal S}_{M}$, $0<\beta\le \alpha$, then by virtue of (\ref{difference_Fourier_Coeff}) and (\ref{Fourier coeff}), for arbitrary numbers
$k\in {\mathbb Z}\setminus\{0\}$ and $h\in [0,\delta]$, we have
$$
    \Big|{[\Delta_h^\alpha f]}\widehat {\ \ }(k)\Big|= 2^\beta \Big|\sin \frac {kh}2\Big|^\beta |1-\mathrm{e}^{-{\mathrm i}kh}|^{\alpha-\beta} |\widehat{f}(k)|
            \le
    \delta^\beta|k|^\beta   |1-\mathrm{e}^{-{\mathrm i}kh}|^{\alpha-\beta}|\widehat{f}(k)|\le \delta^\beta  \Big|{[\Delta_h^{\alpha-\beta} f^{(\beta)}]}\widehat {\ \ }(k)\Big|,
$$
and therefore property (vi) holds.

If ${\alpha\in \mathbb{N}}$ and ${p \in \mathbb{N}},$ then using the representation
$$
    \Delta_{ph}^\alpha f(x)= \sum_{k_1=0}^{p-1} \ldots \sum_{k_\alpha=0}^{p-1}
    \Delta_h^\alpha f(x-(k_1+k_2+\ldots+k_\alpha) h),
$$
and the relation
\begin{eqnarray}\nonumber
    \Big| [\Delta_h^\alpha f(x-(k_1+\ldots+k_\alpha) h)]\widehat { \ \ } (k) \Big|
           & =&
    \Big|\frac{1}{2\pi} \int_{-\pi}^\pi \sum\limits_{j=0}^\alpha (-1)^{j} {\alpha \choose j} f_{jh}(x-(k_1+\ldots+k_\alpha) h) \mathrm{e}^{-\mathrm{i}kx}~{\mathrm d}x\Big|
\\ \nonumber
    &\le& \Big|\frac{1}{2\pi} \int_{-\pi}^\pi \sum\limits_{j=0}^\alpha (-1)^{j} {\alpha \choose j} f_{jh}(x) \mathrm{e}^{-\mathrm{i}kx}~{\mathrm d}x\Big|=\Big| [\Delta_h^\alpha f(x)]\widehat { \ \ } (k) \Big|,
\end{eqnarray}
we get
\begin{eqnarray}\nonumber
    \| \Delta_{ph}^\alpha f(x) \|_{_{\scriptstyle  M}} &=&
    \inf\Big\{a>0: \sum\limits_{k\in \mathbb{Z}} M\Big(  \Big| \sum_{k_1=0}^{p-1}   \ldots \sum_{k_\alpha=0}^{p-1} [\Delta_h^\alpha f(x-(k_1+\ldots+k_\alpha) h)]\widehat { \ \ } (k)  \Big| /a   \Big) \le 1 \Big\}
\\ \nonumber
    &\le& \inf\Big\{a>0: \sum\limits_{k\in \mathbb{Z}}
    M\Big( p^\alpha \Big| [\Delta_h^\alpha f(x)]\widehat { \ \ } (k)  \Big| /a   \Big) \le 1 \Big\}\le
    p^\alpha \| \Delta_{h}^\alpha f(x) \|_M,
\end{eqnarray}
and property (vii) is proved. To prove (viii) it is sufficient to consider the case $\delta< \eta$
(for $\delta \ge \eta$,  property (viii) is obvious). Choosing the number $p$ such that ${\eta \over \delta} \le p < {\eta \over \delta}+1$, by virtue (i) and (vii), we obtain
$$
    \omega_\alpha(f; \eta) \le \omega_\alpha(f; p \delta)_{_{\scriptstyle  M}} \le p^\alpha \omega_\alpha(f;  \delta)_{_{\scriptstyle  M}} \le
    ({\eta \over \delta}+1)^{\alpha}   \omega_\alpha(f, \delta)_{_{\scriptstyle  M}}.
$$
\vskip -3mm$\hfill\Box$

\vskip 2mm \textit{ Proof of Lemma \ref{Lemma_3}.} The right-hand side of inequality (\ref{estim-for-norms})
is obtained from the Young inequality $uv\le M(u)+\tilde{M}(v)$, where $u,v \ge 0$,
(see, e.g. \cite[Ch.~1, \S 2]{Krasnosel'skii_Rutickii_M1961}) as follows
\begin{eqnarray}\nonumber
    \| {f\|^\ast_{_{\scriptstyle  M}}/\|f\|_{_{\scriptstyle  M}}} &=&
    \Big\| {f/\|f\|_{_{\scriptstyle  M}}} \Big\|^\ast_{_{\scriptstyle  M}}=
    \sup \Big\{ \sum_{k \in \mathbb{Z}} { \lambda_k|\widehat{f}(k)}|/
    {\|f\|_{_{\scriptstyle  M}}} : \  \lambda\in \Lambda\Big\}
\\ \nonumber
    &\le& \sup \Big\{\sum\limits_{k \in \mathbb{Z}}
    \Big(  M(|{\widehat{f}(k)|}/{\|f\|_{_{\scriptstyle  M}}})+
    \tilde{M} (\lambda_k)\Big): \ \lambda\in \Lambda\Big\}\le 2.
\end{eqnarray}

To prove the left-hand side of the inequality (\ref{estim-for-norms}), we choose an arbitrary
function $f \in {\mathcal S}_{M}$ such that
$\|f\|^\ast_{_{\scriptstyle  M}}=1$, and show that for this function the inequality
$\|f\|_{_{\scriptstyle  M}}\le 1$ holds.

Using the relation $M(u)=\int_0^{u} p (t)~\mathrm{d}t$,  $u\ge 0$, we define the function $p =p(t)$, $t\ge 0$,
and consider the sequence  $\lambda^*=\{\lambda_k^*\}_{k\in \mathbb{Z}}$, where $\lambda_k^* :=p(|\widehat{f}(k)|)$.
Then for any  $k\in \mathbb{Z}$, the inequality
\begin{equation} \label{Yung-equal}
    \lambda_k^*|\widehat{f}(k)| =M(|\widehat{f}(k)|)+\tilde{M}(\lambda_k^*).
\end{equation}
holds  (see \cite{Krasnosel'skii_Rutickii_M1961}). Also note that for any Orlicz function  $M$ the function $\tilde{M}$, defined by (\ref{M_tilde}), is also convex (see \cite{Krasnosel'skii_Rutickii_M1961}) and  satisfies the inequality
\begin{equation} \label{convex}
    \tilde{M}(\mu v)\le \mu\tilde{M}(v),\quad 0\le \mu\le 1.
\end{equation}
If we assume that $\sum_{k\in \mathbb{Z}} \tilde{M}(\lambda_k^*) >1$, then by  (\ref{convex}) we see that
\begin{equation} \label{relation_le1}
   \sum_{k\in \mathbb{Z}} \tilde{M}\Bigg( {\lambda_k^* \over \sum_{j\in \mathbb{Z}} \tilde{M}(\lambda_j^*)}\Bigg)\le \sum_{k\in \mathbb{Z}}{ \tilde{M} ( \lambda_k^*) \over
    \sum_{j\in \mathbb{Z}} \tilde{M}(\lambda_j^*) } \le 1.
\end{equation}
Taking into account (\ref{Yung-equal}), (\ref{relation_le1}), the definition of the set $\Lambda$ and the equality $\|f\|^\ast_{_{\scriptstyle  M}}=1,$ we get the contradiction
\begin{eqnarray}\nonumber
    \sum\limits_{k\in \mathbb{Z}} M(|\widehat{f}(k)|)+
    \sum\limits_{k\in \mathbb{Z}}\tilde{M}(\lambda_k^*)&=&
    \sum\limits_{k\in \mathbb{Z}}\lambda_k^* |\widehat{f}(k)|=    \sum\limits_{i \in \mathbb{Z}}\tilde{M}(\lambda_i^*)\sum\limits_{k\in \mathbb{Z}}
    |\widehat{f}(k)| \frac{\lambda_k^*}{\sum _{i\in \mathbb{Z}}\tilde{M}(\lambda_i^*)}
\\ \nonumber
    &\le&  \sum\limits_{i \in \mathbb{Z}}\tilde{M}(\lambda_i^*)
    \sup\Big\{ \sum\limits_{k \in \mathbb{Z}}\lambda_k|\widehat{f}(k)| : \quad\lambda\in \Lambda\Big\}=
     \sum\limits_{k\in \mathbb{Z}}\tilde{M}(\lambda_k^*).
\end{eqnarray}
Consequently, $\sum_{k\in \mathbb{Z}} \tilde{M}(\lambda_k^*) \le 1$ and therefore $\lambda_k^*\in \Lambda$. Then, taking into account (\ref {Yung-equal}), we obtain
$$
    \sum\limits_{k\in \mathbb{Z}}M (\widehat{f}(k)) \le
     \sum\limits_{k\in \mathbb{Z}}|\widehat{f}(k)| \lambda_k^*\le \|f\|^\ast_{_{\scriptstyle  M}}=1,
$$
hence, $\|f\|_{_{\scriptstyle  M}}\le 1.$  $\hfill\Box$

\vskip 2mm \textit{ Proof of Lemma \ref{Lemma_4}.} Since for any polynomial of the form $\tau_{n}(x){=}\sum_{|k|\le n}  c_{k}\mathrm{e}^{\mathrm{i}kx}$ we have
$\|\tau_{n}^{(\alpha)}\|_{_{\scriptstyle  M}}=\inf\{a>0: \sum_{|k|\le n}M(|k|^\alpha |c_{k}|/a)\le 1\}$, then similarly to (\ref{Trigon_P}), we obtain
\[
    \sum_{|k|\le n} M\Big(|{[\Delta_h^\alpha \tau_{n}]}\widehat {\ \ }(k)|/a_1\Big)\le
     \sum_{|k|\le n} M\Big( |kh|^\alpha |c_{k}|/a_1\Big) \le \sum_{|k|\le n}M\Big(|k|^\alpha |c_{k}|/\|\tau_{n}^{(\alpha)}\|_{_{\scriptstyle  M}}\Big)\le 1,
\]
when $a_1:=|h|^\alpha \|\tau_{n}^{(\alpha)}\|_{_{\scriptstyle  M}}$. Therefore,  $\|\Delta_h^{\alpha} \tau_{n}\|_{_{\scriptstyle  M}}\le |h|^\alpha \|\tau_{n}^{(\alpha)}\|_{_{\scriptstyle  M}}$.

In (\ref{Bermstain-inequl-gener}), the first inequality  is trivial in the cases where $h=0$ or $|h|=2\pi/n$. So, now let $0<|h|<2\pi/n$. Since
 \[
 \|\Delta_h^{\alpha} \tau_{n}\|_{_{\scriptstyle  M}}=\inf\Big\{a>0: \sum_{|k|\le n} M\Big(2^\alpha \Big|\sin \frac {kh}2\Big|^\alpha |c_{k}|/a\Big)\le 1\Big\}
 \]
and  the function $t/\sin t$ increase on $(0,\pi)$, then for  $a_2:=\Big|\frac{n/2}{\sin (nh/2)}\Big|^\alpha \|\Delta_h^{\alpha} \tau_{n}\|_{_{\scriptstyle  M}}$ we get
 \[
 \sum_{|k|\le n}M(|k|^\alpha |c_{k}|/{a_2})=
  \sum_{|k|\le n}M\Big(\Big|\frac{kh/2}{\sin (kh/2)}\Big|^\alpha
  \Big|\frac{\sin (kh/2)}{h/2}\Big|^\alpha
   |c_{k}|/{a_2}\Big)
 \]
 \[
   \le\sum_{|k|\le n}M\Big(\Big|\frac{nh/2}{\sin (nh/2)}\Big|^\alpha
  \Big|\frac{\sin (kh/2)}{h/2}\Big|^\alpha
   |c_{k}|/{a_2}\Big)=\sum_{|k|\le n} M\Big(2^\alpha \Big|\sin \frac {kh}2\Big|^\alpha |c_{k}|/\|\Delta_h^{\alpha} \tau_{n}\|_{_{\scriptstyle  M}}\Big)\le 1.
 \]
Thus, the first inequality in (\ref{Bermstain-inequl-gener}) also holds.  $\hfill\Box$

\end{document}